\documentclass[12pt]{amsart}
\usepackage{amssymb,amsmath,amsthm}
\usepackage{a4}
\usepackage{color}
\usepackage{enumitem}
\usepackage[colorlinks=true,
linkcolor=webgreen,
filecolor=webbrown,
citecolor=webgreen]{hyperref}
\definecolor{webgreen}{rgb}{0,0,1}
\definecolor{recrown}{rgb}{1,.2,.6}
\usepackage{fullpage}
\usepackage{orcidlink}

\begin{document}
\newtheorem{theorem}{Theorem}
\newtheorem{corollary}[theorem]{Corollary}
\newtheorem{lemma}[theorem]{Lemma}
\theoremstyle{definition}
\newtheorem{example}{\bf Example}
\theoremstyle{theorem}
\newtheorem{conjecture}[theorem]{Conjecture}
\newtheorem{thmx}{\bf Theorem}
\renewcommand{\thethmx}{\text{\Alph{thmx}}}
\newtheorem{lemmax}{\bf Lemma}
\renewcommand{\thelemmax}{\text{\Alph{lemmax}}}
\hoffset=-0cm
\theoremstyle{definition}
\newtheorem*{definition}{Definition}
\theoremstyle{remark}
\newtheorem*{remark}{\bf Remark}
\theoremstyle{remark}
\newtheorem*{remarks}{\bf Remarks}
\parindent=0cm
\title{\bf A generalized Dumas irreducibility criterion}
\author{Rishu Garg$^{1}$ {\large \orcidlink{0009-0008-2348-8340}}}
\author{Jitender Singh$^{2,\dagger}$ {\large \orcidlink{0000-0003-3706-8239}}}
\address[1,2]{Department of Mathematics,
Guru Nanak Dev University, Amritsar-143005, India\newline 
{\tt rishugarg128@gmail.com}, {\tt jitender.math@gndu.ac.in}}

\markright{}
\date{}
\footnotetext[2]{Corresponding author email: {\tt jitender.math@gndu.ac.in}

2020MSC: {Primary 12E05; 11C08}\\

}
\maketitle
\newcommand{\K}{\mathbb{K}}
\newcommand{\dd}{\Gamma^{\text{div}}}
\newcommand{\dv}{G_v^{\text{div}}}
\begin{abstract}
As an extension of the classical irreducibility result of Dumas, a factorization result for polynomials over any valued field with a Krull valuation of arbitrary rank is proved. Further, among other results, a lower degree factor bound on factors of a given polynomial over a valued field with a Krull valuation is proved. These factorization results not only unify several known irreducibility results for polynomials over the said domains but also provide us sharp bounds on degrees of irreducible factors of the underlying polynomials.
\end{abstract}
\section{Introduction}
The problem of determining the irreducibility of polynomials over a given algebraic domain has been fundamental in algebra and number theory, since the nineteenth century. Symmetry of a geometrical or natural object is  measured through abstract groups in mathematics. In view of this, symmetry of roots of an irreducible polynomial is captured by the action of the Galois group associated with the underlying polynomial. Moreover, irreducible polynomials are used in constructing finite fields, which in turn play a key role in data encryption in computer science, which makes such studies important from technological point of view as well.

One of the earliest known classical irreducibility results is due to Sch\"onemann \cite{S} in 1846, who established a criterion for testing irreducibility of certain classes of polynomials having integer coefficients. A few years later, Eisenstein \cite{E} in 1850  formulated his irreducibility criterion, which turns out to be a particular case of Sch\"onemann's irreducibility result. For an historical account of Eisenstien's irreducibility criterion, the reader is referred to see the excellent exposition in \cite{Cox}. In 1906, Dumas  \cite{D} proved a far reaching generalization of the irreducibility criteria of Sch\"onemann and  Eisenstein by using the geometric tool of Newton polygons. Over the years, these exquisite irreducibility criteria have been extended, generalized, and reformulated using Newton polygons and valuation theory. Using Newton polygon approach, Bonciocat \cite{B2015} obtained a number of new irreducibility results for polynomials having integer coefficients which satisfy divisibility properties with regard to two or more distinct prime numbers in an analogy with Dumas irreducibility criterion. Using valuation theory,  Khanduja and Saha \cite{KS} extended Dumas irreducibility criterion for polynomials over valued fields with a valuation of arbitrary rank.  Another generalization of Dumas irreducibility criterion was obtained in \cite[Theorem 1.1]{KK} for polynomials over the Krull valued fields. In 2013, Weintraub \cite{W} gave a mild generalization of Eisenstein's irreducibility criterion, which also provided a correction to a false claim made by Eisenstein himself in \cite{E}. Building upon this, in \cite{JK}, the authors proved Weintraub's result under weaker assumptions and in a more general setting of valued fields, which also extends the irreducibility criterion of Dumas. To state their result, we need some definitions and notations described as follows. Let $\Gamma$ be an ordered abelian group with a total order $\leq$ and binary operation $+$. Let $\infty$ be a symbol such that for all $x\in \Gamma$, $\infty=\infty+\infty=x+\infty=\infty+x$. For a subgroup $G$ of $\Gamma$, a Krull valuation $v$ on a field $\K$ is a surjective map $v:\K\rightarrow G\cup \{\infty\}$ satisfying the following axioms:
\begin{enumerate}
\item[$(i)$] $v(x)=\infty$ if and only if $x=0$,~$(ii)$ $v(xy)=v(x)+v(y)$ for all $x,y\in \K$,
\item[$(iii)$] $v(x+y)\geq \min\{v(x),v(y)\}$ for all $x,y\in\K$.
\end{enumerate}
The group $G=v(\K^\times)$ is called the value group of $v$. In view of the above definition we observe that $v(-x)=v(x)$ and $v(x+y)=\min\{v(x),v(y)\}$ if and only if $v(x)\neq v(y)$.
Using these notions,   the following irreducibility result was proved in \cite{JK}.
\begin{thmx}[Jhorar and Khanduja \cite{JK}]\label{th:A}
Let $v$ be a Krull valuation of arbitrary rank of a field $\K$ with its value group $G_v$. Let $f=a_0+a_1z+\cdots +a_nz^n \in \K[z]$ be a polynomial of degree $n$ and $k$ be a smallest nonnegative integer for which the following conditions are satisfied.
\begin{enumerate}
	\item[$(i)$] $v(a_n)=0$,
\item [$(ii)$] $\frac{v(a_k)}{n-k}<\frac{v(a_i)}{n-i}$ for all $i=0,\ldots, n-1$ with $i\neq k$,
	\item[$(iii)$] $v(a_k)\notin dG_v$ for any positive divisor $d>1$ of $n-k$.
\end{enumerate}	
Then any factorization $f(z)=f_1(z)f_2(z)$ of $f$ in $\K[z]$ has a factor of degree $\leq k$. In particular, if $k=0$, then the polynomial $f$ is irreducible in $\K[z]$.
\end{thmx}
In view of the hypothesis (i) of Theorem \ref{th:A}, if  valuation of a coefficient other than the leading coefficient is zero, then on possibly modifying the hypothesis (ii) can one still have any information about the factorization of $f?$ This question has been answered affirmatively in \cite[Theorem 1]{JRG2} for polynomials over $\mathbb{Z}$ and $v=v_p$, the $p$-adic valuation of $\mathbb{Q}$. Recently, in \cite{Singh25}, the author proved the following factorization result for polynomials over discrete valuation domains using Newton-polygon approach.
\begin{thmx}[Singh \cite{Singh25}]\label{th:B}
Let $(R,v)$  be a discrete valuation domain. Let $f=a_0 + a_1z + \cdots + a_nz^n \in {R}[z]$ be a primitive polynomial such that there exist indices $j$ and $k$ with $1\leq k+1 \leq j\leq n$ for which the following hold.
\begin{enumerate}[label=$(\roman*)$]
\item[$(i)$] $v(a_j)=0$,
\item[$(ii)$]  $\frac{v(a_k)}{j-k}<\frac{v(a_i)}{j-i}$ for all $i=0,1,\ldots, j-1$ with $i\neq k$,
\item[$(iii)$] $\gcd(v(a_k),~j-k) = 1$.
\end{enumerate}
Then any factorization $f(z)=f_1(z)f_2(z)$ of $f$ in $R[z]$ has a factor of degree $\leq n-j+k$. In particular, if $j=n$ and $k=0$, then the polynomial $f$ is irreducible in $R[z]$.
\end{thmx}
The case  $j=n$ of Theorem \ref{th:B} is the irreducibility result of Weintraub \cite{W} for polynomials over discrete valuation domains. As a generalization of Theorems \ref{th:A}-\ref{th:B}, we may have the following  result.
\begin{theorem}\label{th:1}
Let $v$ be a Krull valuation of arbitrary rank of a field $\K$ with the value group $G_v$. Let $f=a_0+a_1 z+\cdots+a_n z^n \in \K[z]$ be a polynomial of degree $n$. Suppose there exist smallest indices $j$ and $k$ with $1\leq k+1\leq j\leq n$ for which the following conditions are satisfied.
\begin{enumerate}[label=$(\roman*)$]
		\item[$(i)$] $v(a_j)=0$,
\item[$(ii)$] $\frac{v(a_k)}{j-k}<\frac{v(a_i)}{j-i}$ for all indices $i$ with $0\leq i\leq j-1$ and $i\neq k$,
\item[(iii)] If $j<n$, then $\frac{v(a_k)}{j-k}>\frac{v(a_i)}{j-i}$ for all indices $i$ with $j+1\leq i\leq n$,
\item[$(iv)$] $v(a_k)\notin dG_v$ for any positive divisor $d>1$ of $j-k$.
	\end{enumerate}
Then any factorization $f(z)=f_1(z)f_2(z)$ of $f$ in $\K[z]$ has a factor of degree $\leq n-j+k$. In particular, if $j=n$ and $k=0$, then the polynomial $f$ is irreducible in $\K[z]$.
\end{theorem}
We note that the case $j=n$ of Theorem \ref{th:1} is precisely Theorem \ref{th:A}, and Theorem \ref{th:B} corresponds to the case when the value group of $v$ is infinite cyclic. As a direct consequence of Theorem \ref{th:1}, we have the following irreducibility result.
\begin{corollary}\label{c:1}
Let $v$ be a valuation of a field $\K$ with its value group $\mathbb{Z}$. Let $f=a_0+a_1z+\cdots +a_nz^n \in \K[z]$ be a polynomial of degree $n$. Suppose there exist smallest indices $j$ and $k$ with $1\leq k+1\leq j\leq n$ for which the following conditions are satisfied.
	\begin{enumerate}[label=$(\roman*)$]
		\item $v(a_j)=0$,
		\item $\frac{v(a_k)}{j-k}<\frac{v(a_i)}{j-i}$ for all indices $i$ with $0\leq i\leq j-1$ and $i\neq k$,
		\item If $j<n$, then $\frac{v(a_k)}{j-k}>\frac{v(a_i)}{j-i}$ for all indices $i$ with $j+1\leq i\leq n$,
		\item $\gcd(v(a_k),j-k)=1$.
	\end{enumerate}	
	Then any factorization $f(z)=f_1(z)f_2(z)$ of $f$ in $\K[z]$ has a factor of degree $\leq n-j+k$. In particular, if $j=n$ and $k=0$, then the polynomial $f$ is irreducible in $\K[z]$.
\end{corollary}
The Corollary \ref{c:1} in particular, yield results of  the papers \cite[Corollary 1.3]{JK} and \cite[Theorem 1]{JRG2} for $j=n$, and for $f\in \mathbb{Z}[z]$, respectively by taking $v=v_p$, the $p$-adic valuation of $\mathbb{Q}$. The case $j=n$ of the  Corollary \ref{c:1} is the irreducibility result of Weintraub \cite{W}. Further, the case $j=n$ and $k=0$ of Corollary \ref{c:1} corresponds to a generalization of the classical irreducibility criterion of Dumas \cite{D}.

Our next result is the following generalization of Weintraub's factorization result \cite{W}.
\begin{theorem}\label{th:2}
Let $v$ be a Krull valuation of arbitrary rank of a field $\K$ with  the value group $G_v$. Let $f=a_0+a_1 z+\cdots+a_nz^n \in \K[z]$ be a polynomial of degree $n$. Suppose there exists a smallest index $j$ with $1\leq j\leq n$ for which the following conditions are satisfied.
\begin{enumerate}[label=$(\roman*)$]
\item[$(i)$] $v(a_j)=0$,
\item[$(ii)$] $\frac{v(a_0)}{j}\leq \frac{v(a_i)}{j-i}$ for all $i=0,1,\ldots, j-1$,
\item[(iii)] If $j<n$, then $\min\Bigl\{\frac{v(a_0)}{j}, \frac{v(a_n)}{j-n}\Bigr\} \geq \frac{v(a_i)}{j-i}$ for all $i=j+1, \ldots, n-1$,
\item[$(iv)$] Let $d_1, d_2$ be the smallest positive integers such that $d_1\frac{v(a_0)}{j}\in G_v$, and if $j<n$, then also $d_2\frac{v(a_n)}{j-n}\in G_v$.
\end{enumerate}
Then any factorization of $f(z)=f_1(z)f_2(z)$ of $f$ in $\K[z]$ has a factor of degree at least $\delta_f$, where  $\delta_f$ is defined as
\begin{eqnarray*}
\delta_f=\begin{cases} \min\{d_1,d_2\}, & \text{if} ~ j < n; \\
d_1, & \text{if} ~ j=n.
\end{cases}
\end{eqnarray*}
\end{theorem}
We may have the following factorization result, which also generalizes one of the main results of the paper \cite[Theorem 3.7]{Kh23}.
\begin{theorem}\label{th:3}
Let $v$ be a Krull valuation of arbitrary rank of a field $\K$ with  the value group $G_v$. Let $f=a_0+a_1 z+\cdots+a_nz^n\in \K[z]$ be a polynomial of degree $n$. Suppose there exists a smallest index $j$ with $1\leq j\leq n$ for which the following conditions are satisfied.
\begin{enumerate}[label=$(\roman*)$]
\item[$(i)$] $v(a_j)=0$,
\item[$(ii)$] $(j-i)v(a_0)\leq jv(a_i)$ for all $i=0,1,\ldots, m-1$, $i\neq j$,
\item[(iii)] If $j<n$, then $(j-i)\frac{v(a_n)}{j-n} \leq v(a_i)$ for all $i=0, \ldots, n-1$, $i\neq j$,
\item[$(iv)$] Let $d_1$ and $d_2$ be the smallest positive integers such that $d_1\frac{v(a_0)}{j}\in G_v$, and if $j<n$, then also $d_2\frac{v(a_n)}{j-n}\in G_v$.
\end{enumerate}
Then each irreducible factor of $f$ in $\K[z]$ has degree at least $\delta_f$, where  $\delta_f$ is defined as
\begin{eqnarray*}
\delta_f=\begin{cases} \min\{d_1,d_2\}, & \text{if} ~ j < n; \\
d_1, & \text{if} ~ j=n.
\end{cases}
\end{eqnarray*}
\end{theorem}
\section{Proofs}
To prove Theorem \ref{th:1}, we essentially use the technique developed in \cite{JK}, and we proceed as follows. Let $(\Gamma,+,\leq)$ be an ordered abelian group. The divisible hull $\dd$ of $\Gamma$ is the set
 \begin{eqnarray*}
 \dd=\Bigl\{\frac{x}{k}~|~x\in \Gamma,~k\in \mathbb{N}\Bigr\}.
 \end{eqnarray*}
 The binary operation $+$ and the order $\leq$ of $\Gamma$ extend to $\dd$ in the usual way so that $(\dd, +, \leq)$ is an ordered abelian group with $\Gamma$ as its ordered abelian subgroup (see \cite{Fuchs1970}).

The following general result will be used in proving Theorem \ref{th:1}.
\begin{lemma}\label{LL1}
    Let $\Gamma$ be an ordered abelian group. Let $g\in \Gamma$ and $N\in \mathbb{N}$ with $N>1$ be such that $g\not\in d \Gamma$ for every positive divisor $d>1$ of $N$. Then $ig\not\in \Gamma$ for each $i=1,\ldots, N-1$.
\end{lemma}
\begin{proof} Assume on the contrary that there exists an index $i$ with $1\leq i\leq N-1$ for which $ig\in \Gamma$. Let $\gamma=g/N\in \dd$. Then $N\gamma=g\in \Gamma$. Let $d=\gcd(i,N)<N$, since $1\leq i<N$. Then we have $(N/d)>1$. By Euclidean algorithm, there exist integers $x$ and $y$ for which $xi+yN=d$. We then have $d\gamma=x(i\gamma)+y(N\gamma)\in \Gamma$. Consequently, we get $g=N\gamma=(N/d)(d\gamma) \in \Gamma$ for the positive divisor $(N/d)>1$ of $N$. This contradicts the hypothesis, and so, the Lemma holds.
\end{proof}
\begin{proof}[\bf Proof of Theorem \ref{th:1}]
Let $\gamma_j=v(a_k)/(j-k)$. Then $\gamma_j\in \dv$ and $(j-k)\gamma_j\in G_v$. This in view of the hypothesis $(iv)$  and Lemma \ref{LL1} tells us that $i\gamma_j\not\in dG_v$ for each divisor $d>1$ of $(j-k)$.

Let $G=\langle G_v,\gamma_j\rangle$ be the subgroup of $\dv$ generated by $G_v$ and $\gamma_j$.
Let $w_j:\K(z)\rightarrow G\cup\{\infty\}$ be such that
\begin{eqnarray*}
w_j\Bigl(\sum_{i}\alpha_i z^i\Bigr) &=& \min_i\{v(\alpha_i)+i\gamma_j\},
\end{eqnarray*}
for all $\sum_{i}\alpha_i z^i\in \K[z]$, and for nonzero $f$ and $g$ in $\K[z]$, we define $w_j(f/g)=w_j(f)-w_j(g)$. Then $w_j$ becomes a valuation of $\K(z)$ with its value group $G$ as follows from \cite[Theorem 2.2]{EP}. By the hypotheses $(ii)$ and $(iii)$, we have for all $i\neq k,~j$ that
\begin{eqnarray*}
v(a_i)+i\gamma_j>\frac{j-i}{j-k}v(a_k)+i\gamma_j=j\gamma_j= v(a_k)+k\gamma_j,
\end{eqnarray*}
which shows that $v(a_k)+k\gamma_j=\min_i\{v(a_i)+i\gamma_j\}=w_j(f)$.

Now suppose on the contrary that $f(z)=f_1(z)f_2(z)$ for nonconstant polynomials $f_1$ and $f_2$ in $\K[z]$ such that $\min\{\deg f_1,\deg f_2\}>n-j+k$. We may write $f_1$ and $f_2$ as
\begin{eqnarray*}
f_1=b_0+b_1z+\cdots+b_{n_1}z^{n_1},~f_2=c_0+c_1z+\cdots+c_{n_2} z^{n_2},~b_{n_1}c_{n_2}\neq 0,~n_1+n_2=n.
\end{eqnarray*}
Since $v(a_j)=0$, there exist smallest indices $r\leq n_1$ and $s\leq n_2$ with $r+s=j$ for which $v(b_rc_s)=0$.
Let $k_i$ denote the smallest index for which $w_j(f_i)$ occurs for each $i=1,2$. We then have
\begin{eqnarray*}
	v(b_i)+i\gamma_j\geq w_j(f_1),~v(c_\ell)+\ell\gamma_j\geq w_j(f_2),
\end{eqnarray*}
for $0\leq i \leq n_1$, $0\leq \ell\leq n_2$ with strict inequality whenever $i<k_1$ or $\ell<k_2$. In view of this, the smallest index for which $w_j(f_1)+w_j(f_2)$ occurs is $k_1+k_2$. Since $w_j(f)=w_j(f_1)+w_j(f_2)$, and the smallest index at which $w_j(f)$ occurs is $k$, we have $k_1+k_2=k$. Since $\min\{n_1,n_2\}>n-j+k$, we deduce that $j-n_i>k$ for each $i=1,2$, and so, we find that
	$k_i\leq k<r+s-n_i,~i=1,2$. These inequalities along with the inequalities $r-n_1\leq 0$ and $s-n_2\leq 0$ yield the following:
\begin{eqnarray}\label{e2a}
	r-k_1>0,~s-k_2>0.
\end{eqnarray}
Since $k_1+k_2=k$, this in view of \eqref{e2a} gives us the following:
\begin{eqnarray*}
    j-k=(r+s)-(k_1+k_2)=(r-k_1)+(s-k_2),
\end{eqnarray*}
which shows that
$0<r-k_1<j-k~\text{and}~0<s-k_2<j-k$.
Since we have
\begin{eqnarray*}
w_j(f_1)&=&v(b_{k_1})+k_1\gamma_j\leq v(b_r)+r \gamma_j,\\
w_j(f_2)&=&v(c_{k_2})+k_2\gamma_j\leq v(c_{s})+s\gamma_j,
\end{eqnarray*}
it follows that $w_j(f_1)+w_j(f_2)\leq v(b_rc_s)+(r+s)\gamma_j=j\gamma_j$. This in view of the fact that $j\gamma_j=w_j(f)$ tells us that
\begin{eqnarray*}
	j\gamma_j=w_j(f)=w_j(f_1f_2)=w_j(f_1)+w_j(f_2)\leq v(b_rc_s)+(r+s)\gamma_j=j\gamma_j,
\end{eqnarray*}
which forces the equalities $w_j(f_1)=v(b_r)+r\gamma_j$ and $w_j(f_2)=v(c_s)+s\gamma_j$. We then have
\begin{eqnarray*}
w_j(f_1)=v(b_r)+r\gamma_j=v(b_{k_1})+k_1\gamma_j,
\end{eqnarray*}
and so, $(r-k_1)\gamma_j =(v(b_{k_1})-v(b_r))\in G_v$, which contradicts Lemma \ref{LL1}.
\end{proof}
\begin{proof}[\bf Proof of Corollary \ref{c:1}]
Observe that the hypothesis $(iv)$ of the corollary implies that for any divisor $d>1$ of $j-k$, $v(a_k)\not\in dG_v=d\mathbb{Z}$, since otherwise, $d$ divides $v(a_k)$, and so, $d$ divides $\gcd(v(a_k),j-k)=1$, which is a contradiction. So, the polynomial $f$ satisfies the hypotheses of Theorem \ref{th:1} for $G_v=\mathbb{Z}$, and the corollary follows.
\end{proof}
\begin{proof}[\bf Proof of Theorem \ref{th:2}]
The case $j=n$ follows from Theorem 1.1 of \cite{AJ}. In view of this, we will prove the theorem for the case when $j<n$. We proceed as follows. Let
\begin{eqnarray*}
\gamma=\min\Bigl\{\frac{v(a_0)}{j},~\frac{v(a_n)}{j-n}\Bigr\}\in \dv.
\end{eqnarray*}
Let	$G=\langle G_v,\gamma\rangle$ be the subgroup of $\dv$ generated by $G_v$ and $\gamma$. Let $w:\K(z)\rightarrow G\cup\{\infty\}$ be such that
\begin{eqnarray*}
w\Bigl(\sum_{i}\alpha_i z^i\Bigr) &=& \min_i\{v(\alpha_i)+i\gamma\},
\end{eqnarray*}
for all $\sum_{i}\alpha_i z^i\in \K[z]$, and for all nonzero $f$ and $g$ in $\K[z]$, we define $w(f/g)=w(f)-w(g)$. Then $w$ becomes a valuation of $\K(z)$ with its value group $G$ as follows from \cite[Theorem 2.2]{EP}. We now have the following two cases.
	
\textbf{Case I}. $\gamma={v(a_0)}/{j}$. In this case, we have
\begin{eqnarray}\label{e7b}
\frac{v(a_0)}{j} \leq \frac{v(a_n)}{j-n}.
\end{eqnarray}
Using \eqref{e7b} along with the hypotheses $(ii)$ and $(iii)$, we have for all $i\neq j$ that
\begin{eqnarray*}
v(a_i)+i\gamma&\geq & \frac{j-i}{j}v(a_0)+i\gamma=j\gamma =v(a_0),
\end{eqnarray*}
which shows that
\begin{eqnarray}\label{e7c}
v(a_0)=\min_i\{v(a_i)+i\gamma\}=w(f).
\end{eqnarray}
Now suppose that $f(z)=f_1(z)f_2(z)$
for nonconstant polynomials $f_1$ and $f_2$ in $\K[z]$. We may write $f_1$ and $f_2$ as
\begin{eqnarray*}
f_{t}=b_{t0}+b_{t1}z+\cdots+b_{tn_t}z^{n_t},~t=1,2,
\end{eqnarray*}
where $n_t\geq 1$, $b_{1n_1}b_{2n_2}\neq 0$ and $n_1+n_2=n$. Since $v(a_j)=0$, there exist smallest index $r_t\leq n_t$ for each $t$ with $r_1+r_2=j$ for which $v(b_{1r_1}b_{2r_2})=0$. Let $k_t$ denote the smallest index for which $w(f_t)$ occurs for each $t=1,2$.  We then have
\begin{eqnarray*}
v(b_{t\ell_t})+\ell_t\gamma\geq w(f_t),
\end{eqnarray*}
for $0\leq \ell_t \leq n_t$, for each $t=1,2$ with strict inequality whenever $\ell_t<k_t$ for  $t=1$, or $t=2$. In view of this, the smallest index for which $w(f_1)+w(f_2)$ occurs is $k_1+k_2$. Since $w(f)=w(f_1)+w(f_2)$, and since from \eqref{e7c}, the smallest index at which $w(f)$ occurs is $0$, we find that $k_1+k_2=0$, which shows that $k_1=k_2=0$. Consequently, we have
\begin{eqnarray}\label{e7d}
w(f_t)=v(b_{t0})\leq v(b_{tr_t})+r_t\gamma,~t=1,2,
\end{eqnarray}
from which we have $w(f_1)+w(f_2)\leq v(b_{1r_1}b_{2r_2})+(r_1+r_2)\gamma=j\gamma$. This in view of the fact that $j\gamma=w(f)$ tells us that
\begin{eqnarray*}
j\gamma=w(f)=w(f_1f_2)= w(f_1)+w(f_2) \leq j\gamma,
\end{eqnarray*}
which in view of \eqref{e7d} forces the equality $w(f_t)=v(b_{t0})=v(b_{tr_t})+r_t\gamma$ for each $t$. Consequently,  we have $r_t\gamma=(v(b_{t0})-v(b_{tr_t})) \in G_v$. Since $r_1+r_2=j>0$, at least one of $r_1$ and $r_2$ is a positive integer. By the minimality of $d_1$, either $n_1\geq r_1\geq d_1$ or $n_2\geq r_2\geq d_1$.
	
\textbf{Case II}. $\gamma={v(a_n)}/{(j-n)}<v(a_0)/j$. In this case, we have
\begin{eqnarray}\label{e7e}
\frac{v(a_n)}{j-n} <\frac{v(a_0)}{j} \leq \frac{v(a_i)}{j-i}, ~\textrm{for}~ 0\leq i\leq j-1.
\end{eqnarray}
Using \eqref{e7e} and the hypotheses $(iii)$, we arrive at the following:
\begin{eqnarray*}
v(a_i)+i\gamma&> & \frac{j-i}{j-n}v(a_n)+i\gamma=j\gamma,
\end{eqnarray*}
for each $0\leq i\leq n-1$, $i\neq j$. So, we have $v(a_i)+i\gamma\geq j\gamma$ for all $i$ with $0\leq i<n$ and $i\neq j$ with the strict inequality if $i<j$. Consequently, we have $v(a_n)+n\gamma =j\gamma=\min_i\{v(a_i)+i\gamma\}=w(f)$, and so, in this case $j$ is the smallest index at which $w(f)$ is attained. Now as before, let $f(z)=f_1(z)f_2(z)$ for nonconstant polynomials $f_1$ and $f_2$ in $\K[z]$, and let
\begin{eqnarray*}
f_{t}=b_{t0}+b_{t1}z+\cdots+b_{tn_t}z^{n_t},
\end{eqnarray*}
where $n_t\geq 1$ for each $t=1,2$, $b_{n_1}c_{n_2}\neq 0$ and $n_1+n_2=n$. Since $v(a_j)=0$, there exist smallest index $r_t\leq n_t$ for each $t$ with $r_1+r_2=j$ for which $v(b_{1r_1}b_{2r_2})=0$. Let $k_t$ denote the smallest index for which $w(f_t)$ occurs for each $t=1,2$. We then have
\begin{eqnarray*}
v(b_{t\ell_t})+\ell_t\gamma\geq w(f_t),
\end{eqnarray*}
for $0\leq \ell_t \leq n_t$, for each $t$ with strict inequality whenever $\ell_t<k_t$ for $t=1$, or $t=2$. In view of this, the smallest index for which $w(f_1)+w(f_2)$ occurs is $k_1+k_2$. This in view of the fact that the smallest index at which $w(f)$ occurs is $j$ tells us that  $k_1+k_2=j$. Since we have
\begin{eqnarray}\label{e7f}
w(f_t)=v(b_{tk_t})+k_t\gamma\leq v(tb_{n_t})+n_t\gamma,~t=1,2,
\end{eqnarray}
it follows that $w(f_1)+w(f_2)\leq v(b_{1n_1}b_{2n_2})+(n_1+n_2)\gamma=v(a_n)+n\gamma$. We then have
\begin{eqnarray*}
v(a_n)+n\gamma=w(f)=w(f_1f_2)=w(f_1)+w(f_2)\leq v(a_n)+n\gamma,
\end{eqnarray*}
which in view of \eqref{e7f} forces the equality $w(f_t)=v(b_{tk_t})+k_t\gamma =v(b_{tn_t})+n_t\gamma$, and so, we have $(n_t-k_t)\gamma=(v(b_{tk_t})-v(b_{tn_t})) \in G_v$ for each such $t$. Since $k_1+k_2=j<n=n_1+n_2$, $n_t-k_t$ is a positive integer for at least one value of $t$. For such an index $t$, we have by the minimality of $d_2$ that $(n_t-k_t)\geq d_2$, and hence $n_t\geq d_2 \geq \delta_f$.
\end{proof}
\begin{proof}[\bf Proof of Theorem \ref{th:3}]
The case $j=n$ is precisely Theorem 1.1 of \cite{AJ}. In view of this, we will prove the theorem for the case when $j<n$. We proceed as follows. Let $\gamma_1= v(a_0)/j\in \dv$ and $\gamma_2= v(a_n)/(j-n)\in \dv$. Let $G_q=\langle G_v,\gamma_q\rangle$ be the subgroup of $\dv$ generated by $G_v$ and $\gamma_q$ for each $q=1,2$. Let $w_q:\K(z)\rightarrow G_q\cup\{\infty\}$ be such that
\begin{eqnarray*}
w_q\Bigl(\sum_{i}\alpha_i z^i\Bigr) &=& \min_i\{v(\alpha_i)+i\gamma_q\},
\end{eqnarray*}
for all $\sum_{i}\alpha_i z^i\in \K[z]$, $q=1,2$, and for nonzero $f$ and $g$ in $\K[z]$, we define $w_q(f/g)=w_q(f)-w_q(g)$. Then $w_q$ becomes a valuation of $\K(z)$ with its value group $G_q$ for each $q$ as follows from \cite[Theorem 2.2]{EP}. By the hypotheses $(ii)$, we have for all $i\neq j$ that
\begin{eqnarray*}
v(a_i)+i\gamma_1&\geq & \frac{j-i}{j}v(a_0)+i\gamma_1=j\gamma_1 =v(a_0),
\end{eqnarray*}
which shows that
\begin{eqnarray}\label{e7g}
v(a_0)=\min_i\{v(a_i)+i\gamma_1\}=w_1(f).
\end{eqnarray}
Now suppose that
\begin{eqnarray*}
f(z)=f_1(z)f_2(z)\cdots f_s(z)
\end{eqnarray*}
be a factorization of $f$ into irreducible factors $f_1,\ldots, f_s$ in $\K[z]$. We may write $f_u(z)=\sum_{t=0}^{n_u}b_{ut}z^t$ for each $u=1,\ldots,s$. Since $v(a_j)=0$, there exist smallest indices $r_u\leq n_u$ for $u=1,2,\ldots, s$ with $r_1+\cdots +r_s=j$ for which $\sum_{u=1}^sv(b_{ur_u})=0$. Let $k_u$ denote the smallest index for which $w_1(f_u)$ occurs for each $u=1,2,\ldots, s$. We then have
$v(b_{ut})+t\gamma_1\geq w_1(f_u)$,
for $0\leq t \leq n_u$, where the  strict inequality holds for $t<k_u$. In view of this, the smallest index for which $w_1(f_1)+w_1(f_2)+\cdots +w_1(f_s)$ occurs is $k_1+k_2+\cdots +k_s$. Since $w_1(f)=w_1(f_1\cdots f_s)=w_1(f_1)+w_1(f_2)+\cdots +w_1(f_s)$, and since from \eqref{e7g}, the smallest index at which $w_1(f)$ occurs is $0$, we find that $k_1+k_2+\cdots +k_s=0$, which shows that $k_1=k_2=\ldots=k_s=0$. Consequently, we have
\begin{eqnarray}\label{e7h}
w_1(f_u)=v(b_{u0})\leq v(b_{ur_u})+r_u\gamma_1,~1\leq u\leq s,
\end{eqnarray}
it follows that $w_1(f_1)+w_1(f_2)+\cdots +w_1(f_s)\leq \sum_{u=1}^s\{v(b_{ur_u})+r_u\gamma_1\}=j\gamma_1$. This in view of the fact that $j\gamma_1=w_1(f)$ tells us that
\begin{eqnarray*}
j\gamma_1=w_1(f)=w_1(f_1f_2\cdots f_s)&=&w_1(f_1)+w_1(f_2)+\cdots+w_1(f_s)=j\gamma_1, 
\end{eqnarray*}
which in view of \eqref{e7h} forces the equalities $w_1(f_u)=v(b_{u0})=v(b_{ur_u})+r_u\gamma_1$ for each $u=1,\ldots,s$, and so, $r_u\gamma_1=(v(b_{u0})-v(b_{ur_u})) \in G_v$. Since $d_1$ is the smallest positive integer such that $d_1\gamma_1 \in G_v$, we deduce that $n_u\geq r_u\geq d_1\geq \delta_f$ for those indices of $u$ for which $r_u>0$.

On the other hand if $r_\ell=0$ for some index $\ell\in \{1,\ldots, s\}$, then we  will show that  $n_\ell\geq d_2$. We proceed as follows. Using the hypotheses $(iii)$, we have \begin{eqnarray*}
v(a_i)+i\gamma_2&> & \frac{j-i}{j-n}v(a_n)+i\gamma_2=j\gamma_2,
\end{eqnarray*}
for all indices $i$ with $0\leq i\leq n-1$ and  $i\neq j$. Consequently, we have $v(a_i)+i\gamma_2\geq j\gamma_2$ for each such $i$ with the strict inequality if $i<j$. We then have $v(a_n)+n\gamma_2 =j\gamma_2=\min_i\{v(a_i)+i\gamma_2\}=w_2(f)$, which shows that $j$ is the smallest index at which $w_2(f)$ is attained.
Let $k'_u$ denote the smallest index for which $w_2(f_u)$ occurs for each $u=1, \ldots, s$. We then have
$v(b_{ut})+t\gamma_2\geq w_2(f_u)$,
for each index $t$ with $0\leq t \leq n_u$ with the strict inequality whenever $t<k'_u$. In view of this, the smallest index for which $w_2(f)= w_2(f_1)+w_2(f_2)+\cdots +w_2(f_s)$ occurs is $k'_1+k'_2+\cdots +k'_s$. Since the smallest index at which $w_2(f)$ occurs is $j$, we must have  $k'_1+k'_2+\cdots +k'_s=j$. Since for each $u=1, \ldots, s$, we have
\begin{eqnarray}\label{e7ia}
w_2(f_u)&=&v(b_{uk'_u})+k'_u\gamma_2\leq v(b_{ur_u})+r_u\gamma_2,\\
\label{e7ib}w_2(f_u)&=&v(b_{uk'_u})+k'_u\gamma_2\leq v(b_{un_u})+n_u\gamma_2.
\end{eqnarray}
 Consequently, we have $w_2(f_1)+w_2(f_2)+\cdots+w_2(f_s)\leq \sum_{u=1}^s\{v(b_{un_u})+n_u\gamma_2\}=v(a_n)+n\gamma_2=j\gamma_2$. This in view of the fact that $v(a_n)+n\gamma_2=j\gamma_2=w_2(f)$ tells us that
\begin{eqnarray*}
j\gamma_2=v(a_n)+n\gamma_2=w_2(f)&=&w_2(f_1f_2\cdots f_s)\\&=& w_2(f_1)+w_2(f_2)+\cdots+w_2(f_s)\\&\leq& v(a_n)+n\gamma_2=j\gamma_2,
\end{eqnarray*}
which in view of \eqref{e7ia}-\eqref{e7ib} forces the equalities $w_2(f_u)=v(b_{uk'_u})+k'_u\gamma_2=v(b_{ur_u})+r_u\gamma_2=v(b_{un_u})+n_u\gamma_2$ for each $u=1,\ldots, s$. We then have $r_u\geq k'_u$ for each $u=1,\ldots, s$. Since $\sum_{u=1}^{s}(r_u-k'_u) =0$ with $r_u\geq k'_u$ for each $u=1,\ldots, s$, we deduce that $r_u=k'_u$ for each $u=1, \ldots, s$. Therefore, $(n_u-r_u)\gamma_2=(v(b_{ur_u})-v(b_{un_u})) \in G_v$ for each such index $u$. In particular, $n_\ell\gamma_2=(n_\ell-r_\ell)\gamma_2\in G_v$, where $n_\ell\geq 1$. This in view of the hypothesis on $d_2$
tells us that $n_\ell\geq d_2\geq \delta_f$.
\end{proof}
\section{Examples}
We now illustrate explicit examples of polynomials whose factorization properties may be deduced using Theorems \ref{th:1} \& \ref{th:2}.

In each of the  examples given below, we will assume that $G_v=\mathbb{Z}\times \mathbb{Z}$, the direct product of the additive group of integers equipped with the dictionary order $\ll$, that is, for  $x,y,z,t\in \mathbb{Z}$, $(x,y)\ll(z,t)$ if either $x<z$, or if $x=z$, then $y<t$.
	\begin{example}
Let $\K=\mathbb{Q}(x)$, the rational function field over $\mathbb{Q}$, and let $p$ be a prime number. For $f=\sum_{i=0}^{n}a_ix^i \in \K$, define $v:\K \rightarrow \mathbb{Z} \times \mathbb{Z}$ by
		\begin{eqnarray*}
			v(f):=\begin{cases} \infty, & \text{if}~ f=0;\\
				(v_p(f),v_\infty(\bar{f})), & \text{if}~ f\neq 0,
			\end{cases}
		\end{eqnarray*}
		where $v_p(f)=\min\{v_p(a_i) ~|~ a_i\neq 0\}$ and $v_\infty$ is the degree valuation on the residue class field $\mathbb{F}_p(x)$. For $f,g \in \K[z]$ with $g\neq 0$, define $v(f/g)=v(f)-v(g)$. Then $v$ is a rank-2 valuation with value group $\mathbb{Z} \times \mathbb{Z}$. Now the polynomial
		\begin{align*}
			P_1 &= (1+4x^4)x+4xz+(1+4x^4)2xz^2+8xz^3+(1+4x^4)2x^2z^4+(1+8x^2+4x^4)z^5+4z^6
		\end{align*}
in $\K[z]$ satisfies the hypothesis of Theorem \ref{th:1} for $n=6$, $j=5$, $k=0$, and $p=2$, since $v(a_5)=(0,0)$, $v(a_0)/(j-0)=(0,-1)/5=(0,-1/5)\ll v(a_i)/(j-i)$ for each $i$ with $1\leq i\leq 4$, $v(a_0)/(j-0)=(0,-1/5)\gg v(a_6)/(j-6)=(-2,0)$, and $v(a_0)=(0,-1)\notin dG_v$ for any positive divisor $d>1$ of $j-k=5$. So, any factorization of $P_1$ in $\K[z]$ has a factor of degree $\leq n-j+k=6-5+0=1$, where we find that $P_1$ factors as
		\begin{eqnarray*}
			P_1=(1+4x^4+4z)(x+2xz^2+2x^2z^4+z^5).
		\end{eqnarray*}
This shows that $P_1$ has the factor $1+4x^4+4z$ of degree 1 in $\K[z]$. Note that Theorem \ref{th:1} applies to $P_1$ whereas  Theorems \ref{th:A} and \ref{th:B} are both inconclusive.
	\end{example}
Let $\K=\mathbb{F}(x,y)$ be the rational function field in two variables over a field $\mathbb{F}$.
For $f=\sum_{t,s}c_{t,s}x^ty^s \in \K$, define $v:\K \rightarrow \mathbb{Z} \times \mathbb{Z}$ by
		\begin{eqnarray*}
			v(f):=\begin{cases} \infty, & \text{if}~ f=0;\\
				\min\{(t,s)~|~c_{t,s}\neq 0\}, & \text{if}~ f\neq 0.
			\end{cases}
		\end{eqnarray*}
Further, for $f,g \in \K[z]$ with $g\neq 0$, define $v(f/g)=v(f)-v(g)$. We observe that $v$ is a rank-2 valuation of $\K$ with its value group $G_v=\mathbb{Z} \times \mathbb{Z}$.
\begin{example} Consider the polynomial
		\begin{align*}
			P_2 &= xy^2-(1-x^2)yz-(1-y-xy)xz^2-(1-x-xy^2)xz^3+(x-y+y^2)xz^4\\
			&\quad- (1-x-x^2)yz^5-(1-xy)z^6+xz^7\in \K[z],
		\end{align*}
which satisfies the hypothesis of Theorem \ref{th:1} for $n=7$, $j=6$, and $k=1$, since $v(a_6)=v(1-xy)=(0,0)$, $v(a_1)/(j-1)=v(-(1-x^2)y)/5=(0,1/5)\ll v(a_i)/(j-i)$ for each $i$ with $0\leq i\leq 5$, $i\neq 1$, $v(a_1)/(j-1)=(0,1/5)\gg  v(a_7)/(j-7)=(-1,0)$, and $v(a_1)=(0,1)\notin dG_v$ for any positive divisor $d>1$ of $j-k=5$. So, any factorization of $f$ in $\K[z]$ has a factor of degree $\leq n-j+k=7-6+1=2$, where we find that $P_2$ factors as
		\begin{eqnarray*}
			P_2=(xy-z+xz^2)(y+xz+xz^2+xyz^3+yz^4+z^5),
		\end{eqnarray*}
which shows that the polynomial $P_2$ has the factor $xy-z+xz^2$ of degree 2 in $\K[z]$. Note that here also, Theorems \ref{th:A} and \ref{th:B} are not applicable in deducing  any information about the factorization of the polynomial $P_2$.
	\end{example}
\begin{example}
    Now consider the polynomial
	\begin{eqnarray*}
		P_3=y+xz+(1+xy^2)z^2+x^2yz^3+xyz^4,
	\end{eqnarray*}
	which satisfies the hypothesis of Theorem \ref{th:2} for $n=4$, $j=2$, and $d_1=d_2=2$ since $v(a_2)=(0,0)$, $v(a_0)/j=(0,1)/2=(0,1/2)\ll v(a_1)/(j-1)=(1,0)$, $\min\{v(a_0)/j,v(a_4)/(j-4)\}=v(a_4)/(j-4)=(-1/2,-1/2)$. $v(a_3)/(j-3)=(-2,-1)\ll(-1/2,-1/2)=v(a_4)/(j-4)$. Further,  
$1\times \frac{v(a_0)}{j}=(0,1/2) \notin G_v$, $1\times \frac{v(a_4)}{j-4}=(-1/2,-1/2)\notin G_v$, whereas $2\times \frac{v(a_0)}{j}=(0,1) \in G_v$ and $2\times \frac{v(a_4)}{j-4}=(-1,-1)\in G_v$. Consequently, any factorization of $f$ in $\K[z]$ has a factor of degree at least $\min\{2,2\}=2$, where we find that the polynomial $P_3$ factors as  
	\begin{eqnarray*}
		P_3=y+xz+(1+xy^2)z^2+x^2yz^3+xyz^4 =(1+xyz^2)(y+xz+z^2).
	\end{eqnarray*}
\end{example}
\subsection*{Disclosure statement}
The author reports to have no competing interests to declare.

\end{document}